\documentclass{amsproc}

\usepackage{euscript}
\usepackage{epsfig}

\newcommand{\CC}{{\Bbb{C}}}

\newcommand{\RR}{{\Bbb{R}}}
\newcommand{\ZZ}{{\Bbb{Z}}}

\newcommand{\calS}{{\EuScript S}}
\newcommand{\eps}{\varepsilon}

\newcommand{\del}{\hat{\delta}}
\newcommand{\delp}{\hat{\delta}'}
\newcommand{\Var}{{\rm Var}^{-1}}

\newtheorem{theorem}{Theorem}[section]
\newtheorem{proposition}[theorem]{Proposition}

\newtheorem{corollary}[theorem]{Corollary}

\theoremstyle{definition}
\newtheorem{definition}[theorem]{Definition}
\newtheorem{example}[theorem]{Example}

\theoremstyle{remark}
\newtheorem{remark}[theorem]{Remark}

\numberwithin{equation}{section}

\begin{document}

\title{On the index of a
vector field at an isolated
singularity}

\author{W.Ebeling}
\address{Institut f\"{u}r Mathematik \\
Universit\"{a}t Hannover \\
Postfach 6009 \\
D-30060 Hannover, Germany\\
e-mail:
ebeling\symbol{'100}math.uni-hannover.de}

\author{S.M.Gusein-Zade}
\address{Deptartment of Mathematics and Mechanics\\
Moscow State University\\
Moscow, 119899, Russia\\
e-mail: sabir\symbol{'100}ium.ips.ras.ru}
\thanks{Partially supported by INTAS--96--0713, DFG 436 RUS 17/147/96, and  RFBR
96--15--96043.}

\subjclass{Primary 57R25, 32S30}
\date{September 9, 1997}

\dedicatory{Dedicated to V.I.Arnold.}

\begin{abstract}
 We consider manifolds with isolated singularities, i.e., topological
spa\-ces which are manifolds (say, $C^\infty$--) outside discrete
subsets (sets of singular points). For (germs of) manifolds with, so called,
cone--like singularities, a notion of the index of an isolated singular
point of a vector field is introduced. There is given a formula for the
index of a gradient vector field on a (real) isolated complete
intersection singularity. The formula is in terms of signatures of
certain quadratic forms on the corresponding spaces of thimbles.
\end{abstract}
\maketitle

\section*{Introduction}\label{sec0}

An isolated singular point of a vector field on $\RR^n$ or on an
$n$-dimensional smooth manifold has a natural integer invariant~---
the index. The formula of Eisenbud, Levine and Khimshiashvili (\cite{EL},
\cite{Kh}) expresses the index of an (algebraically) isolated
singular point of a vector field as the signature of a quadratic form
on a local algebra associated with the singular point. For a
singular point of a gradient vector field there is a formula in terms
of signatures of certain quadratic forms defined by the action of the
complex conjugation on the corresponding Milnor lattice (\cite{GZ},
\cite{V}).

We define a generalisation of the notion of the index of an isolated
singular point of a vector field on a manifold with
isolated cone-like singularities in such a way that the Poincar\'e--Hopf
theorem (the sum of indices of singular points of a vector field on a
closed manifold is equal to its Euler characteristic) holds. It seems
that this notion (though very natural) cannot be found in the
literature in an explicit form. In particular, the index is defined for
vector fields on a germ of a real algebraic variety with an isolated
singularity. We give a generalisation of a formula from
\cite{GZ} for the gradient vector field on an isolated complete
intersection singularity. For that we define (in a somewhat formal way)
the notion of the variation operator for a complete intersection
singularity. We show that this operator is an invariant of the
singularity.

\section{Basic definitions}\label{sec1}

A {\em manifold with isolated singularities} is a topological
space $M$ which has the structure of a smooth (say, $C^\infty$--) manifold
outside of a discrete set $S$ (the {\em set of singular points} of $M$).
A {\em diffeomorphism} between two such manifolds is a homeomorphism
which sends the set of singular points onto the set of singular points and
is a diffeomorphism outside of them. We say that $M$ has a {\em
cone-like singularity} at a (singular) point $P\in S$ if there exists a
neighbourhood of the point $P$ diffeomorphic to the cone over a smooth
manifold $W_P$ ($W_P$ is called the {\em link} of the point $P$). In
what follows we assume all manifolds to have only cone-like
singularities. A (smooth or continuous) {\em vector field} on a
manifold $M$ with isolated singularities is a (smooth or continuous) vector
field on the set $M\setminus S$ of regular points of $M$. The {\em set
of singular points} $S_X$ of a vector field $X$ on a (singular) manifold
$M$ is the union of the set of usual singular points of $X$ on $M\setminus
S$ (i.e., points at which $X$ tends to zero) and of the set $S$ of singular
points of $M$ itself.

For an isolated {\em usual} singular point $P$ of a vector field $X$
there is defined its index $\mbox{ind}_PX$ (the degree of the map
$X/\Vert X\Vert:\partial B\to S^{n-1}$ of the boundary of a small ball $B$
centred at the point $P$ in a coordinate neighbourhood of $P$;
$n=\mbox{dim}\,M$). If the manifold $M$ is closed and has no singularities
($S=\emptyset$) and the vector field $X$ on $M$ has only isolated
singularities, then
\begin{equation}\label{eq1}
\sum_{P\in S_X}\mbox{ind}_PX=\chi(M)
\end{equation}
($\chi(M)$ is the Euler characteristic of $M$).

Let $(M, P)$ be a cone-like singularity (i.e., a germ of a manifold
with such a singular point) and let $X$ be a vector field defined on an
open neighbourhood $U$ of the point $P$. Suppose that $X$ has no singular
points on $U\setminus\{P\}$. Let $V$ be a closed cone--like neighbourhood of
$P$ in $U$ ($V\cong CW_P$, $V\subset U$). On the cone $CW_P=(I\times
W_P)/(\{0\}\times W_P)$ ($I=[0, 1]$) there is defined a natural vector field
$\partial/\partial t$ ($t$ is the coordinate on $I$). Let $X_{rad}$ be the
corresponding vector field on $V$. Let $\widetilde X$ be a smooth vector
field on $U$ which coincides with $X$ near the boundary $\partial U$ of
the neighbourhood $U$ and with $X_{rad}$ on $V$ and has only isolated
singular points.

\begin{definition} The {\em index} $\mbox{ind}_PX$ of the vector field $X$
at the point $P$ is equal to
$$1+\sum_{Q\in S_{\widetilde X}\setminus\{P\}}\mbox{ind}_Q\widetilde X$$
(the sum is over all singular points $Q$ of $\widetilde X$ except $P$ itself).
\end{definition}

For a cone-like singularity at a point $P\in S$, the link $W_P$ and thus
the cone structure of a neighbourhood are, generally speaking, not
well-defined (cones over different manifolds may be {\em locally}
diffeomorphic). However it is not difficult to show
that the index $\mbox{ind}_PX$ does not depend on the choice of a cone
structure on a neighbourhood and on the choice of the vector field
$\widetilde X$.

\begin{example} The index of the ``radial" vector field $X_{rad}$ is equal
to $1$. The index of the vector field $(-X_{rad})$ is equal to $1-\chi(W_P)$
where $W_P$ is the link of the singular point $P$.
\end{example}

\begin{proposition}\label{prop1}
For a vector field $X$ with isolated singular points on a closed manifold $M$
with isolated singularities, the relation {\rm(\ref{eq1})} holds.
\end{proposition}

\begin{definition} One says that a singular point $P$ of a manifold $M$
(locally diffeomorphic to the cone $CW_P$ over a manifold $W_P$) is
{\em smoothable} if $W_P$ is the boundary of a smooth compact manifold.
\end{definition}

The class of smoothable singularities includes, in particular, the class of
(real) isolated complete intersection singularities. For such a singularity,
there is a distinguished cone--like structure on its neighbourhood.

Let $(M, P)$ be a smoothable singularity (i.e., a germ of a manifold
with such a singular point) and let $X$ be a vector field on $(M, P)$
with an isolated singular point at $P$. Let $V=CW_P$ be a closed
cone--like neighbourhood of the point $P$; $X$ is supposed to have no singular
points on $V\setminus\{P\}$. Let the link $W_P$ of the point $P$ be the
boundary of a compact manifold $\widetilde V_P$. Using a smoothing one can
consider the union
$\widetilde V_P\cup_{W_P}(W_P\times[1/2, 1])$ of $\widetilde V_P$ and
$W_P\times[1/2, 1]\subset CW_P$ with the natural identification of
$\partial\widetilde V_P=W_P$ with $W_P\times\{1/2\}$ as a smooth manifold
(with the boundary $W_P\times\{1\}$). The restriction of the vector field
$X$ to $W_P\times[1/2, 1]\subset CW_P$ can be extended to a smooth vector
field $\widetilde X$ on $\widetilde V_P\cup_{W_P}(W_P\times[1/2, 1])$ with
isolated singular points.

\begin{proposition}\label{prop2}
The index $\mbox{ind}_PX$ of the vector field $X$ at the
point $P$ is equal to
$$\sum_{Q\in S_{\widetilde X}}\mbox{ind}_Q\widetilde X -
\chi(\widetilde V_P)+1$$
(the sum is over all singular points of $\widetilde X$ on $\widetilde
V_P$).
\end{proposition}

\begin{remark} In \cite{S}, \cite{GSV} there was defined a notion of the
index of a vector field at an isolated singular point of a complex variety
(satisfying some conditions). That definition does not coincide with the
one given here. These definitions differ by the Euler characteristic of the
smoothing of the singularity of the variety. One can say that the
index of \cite{S}, \cite{GSV} depends on the Euler characteristic of
a smoothing and thus is well-defined only for a singularity with
well-defined topological type of a smoothing (at least with well-defined
Euler characteristic of it). It is valid, e.g., for {\em complex}
isolated complete intersection singularities.

A closely related notion has been discussed in \cite{BG}. That
notion can be considered as a relative version of the index defined here.
After a
previous version of this paper had been submitted and put on the Duke
preprint server
as alg-geom/9710008, the authors' attention was drawn to the preprint
\cite{ASV},
where a somewhat more general notion is defined, which coincides with the index
considered here for real analytic varieties with isolated singularities.
\end{remark}

A generic (smooth or continuous) vector field
on a (singular) analytic variety has zeroes only at isolated points. Thus
it is desirable to have a definition of the index of such
a point. One can use the following definition.

Let $(V,0)\subset (\RR^N, 0)$ be a germ of a real algebraic variety
and let $X$ be a continuous vector field on $(V,0)$ (i.e., the
restriction of a continuous vector field on $(\RR^N, 0)$ tangent to
$V$ at each point) which has an isolated zero at the origin (in $V$).
Let $\calS=\{\Xi\}$ be a semianalytic Whitney stratification of $V$
such that its only zero-dimensional stratum $\Xi^0$ consists of the origin.
Let $\Xi$ be a stratum of the stratification $\calS$ and let $Q$ be a point
of $\Xi$.
A neighbourhood of the point $Q$ in $V$ is diffeomorphic to the direct
product of a linear space $\RR^k$ (the dimension $k$ of which is equal to the
dimension of the stratum $\Xi$) and the cone $CW_Q$ over a compact
singular analytic variety $W_Q$. (A diffeomorphism between two
stratified spaces is a homeomorphism which is a diffeomorphism on
each stratum.) In particular a neighbourhood $U(0)$ of the origin is
diffeomorphic to the cone $CW_0$ over a singular variety $W_0$.
It is not difficult to show that there exists a
(continuous) vector field $\widetilde X$ on $(V, 0)$
such that:
\begin{enumerate}
\item the vector field $\widetilde X$ is defined on the neighbourhood
$U(0)\cong CW_0$ of the origin;
\item $\widetilde X$ coincides with the vector field $X$ in a
neighbourhood of the base $\{1\}\times W_0$ of the cone $CW_0$;
\item the vector field $\widetilde X$ has only a finite number of zeroes;
\item each point $Q\in U(0)$ with $\widetilde X(Q)=0$ has a neighbourhood
diffeomorphic to $(\RR^k,0)\times CW_Q$ in which $\widetilde X(y,z)$
($y\in \RR^k$, $z\in CW_q$)
is of the form $Y(y)+Z_{rad}(z)$, where $Y$ is a germ of a vector
field on $(\RR^k,0)$ with an isolated singular point at the origin,
$Z_{rad}$ is the radial vector field on the cone $CW_Q$.
\end{enumerate}

Let $S_{\widetilde X}$ be the set of zeroes of the vector field
$\widetilde X$ (including the origin). For a point $Q\in S_{\widetilde X}$,
let $\widetilde{ind}(Q):=ind_0 Y$, where $Y$ is the vector field on $(\RR^k, 0)$
described above. We define $ind(0)$ to be equal to $1$ (in this case
$k=0$).

\begin{definition}
$ind_{(V,0)}X=\sum\limits_{Q\in S_{\widetilde X}}
\widetilde{ind}(Q)$.
\end{definition}

\section{On the topology of isolated complete intersection
singularities}\label{sec2}

Let $(V,0) \subset (\CC^{n+p},0)$ be an $(n-1)$-dimensional isolated complete
intersection singularity (abbreviated {\em icis} in the sequel) defined by a
germ of an analytic mapping
$$f=(f_1, \ldots , f_{p+1}): (\CC^{n+p},0) \to (\CC^{p+1},0).$$
(We use somewhat strange notations for the dimension and the number of
equations in order to be consistent with the notations in
Section~\ref{sec3}.) For $\delta > 0$, let
$B_\delta$ be the ball of radius
$\delta$ around the origin in $\CC^{n+p}$. For $\delta > 0$ small enough and
for a generic $t \in
\CC^{p+1}$ with $0< \| t \| << \delta$, the set $$V_t = f^{-1}(t) \cap
B_\delta$$ is a manifold with boundary and is called a {\em Milnor fibre} of
the {\em icis}
$(V,0)$ (or of the germ $f$). The diffeomorphism type of $V_t$ does not depend
on $t$. The manifold $V_t$ is homotopy equivalent to the bouquet of $\mu$
spheres of dimension $(n-1)$, where $\mu$ is the Milnor number of the
{\em icis}
$(V,0)$.

For $t \in \CC^{p+1}$, $0 \leq i \leq p$ we define
\begin{eqnarray*}
(V^{(i)},0) & := & (\{ x \in B_\delta: f_1(x)= \ldots = f_{p-i+1}(x) =
0 \},0), \\
V^{(i)}_t & := & \{ x \in B_\delta: f_j(x)=t_j, 1 \leq j \leq p-i+1 \}
\end{eqnarray*}
and we set 
$(V^{(p+1)}, 0) := (\CC^{n+p}, 0)$.
We assume that $(f_1, \ldots , f_{p+1})$
is a system of functions such that for $0\leq i \leq p $ the germ
$(V^{(i)},0)$ is an $(n+i-1)$-dimensional {\em icis}. For any $t \in
\CC^{p+1}$ with $0 < | t_1 | << |t_2| << \ldots << |t_{p+1} | << \delta$, the
set $V^{(i)}_t$ is the Milnor fibre of the {\em icis} $(V^{(i)},0)$. Here the
condition $0 < | t_1 | << |t_2| << \ldots << |t_{p+1} |$ means that
$t_1, \ldots ,t_{p+1}$ have to be chosen in such a way that for each $i$,
$1 \leq i \leq p+1$, all the critical values of the function $f_i$ on
$V^{(p-i+2)}_t$ are contained in the disc of radius $|t_i|$ around $0$.
We put
\begin{eqnarray*}
\hat{H}^{(i)} & := & H_{n+i}(V^{(i+1)}_t,V^{(i)}_t)
\quad \mbox{for} \  0 \leq i \leq p-1 , \\
\hat{H}^{(p)}  & := & H_{n+p}(B_\delta, V^{(p)}_t).
\end{eqnarray*}

We have short exact sequences (cf.\ \cite{Eb}, \cite{AGV}):
$$
\begin{array}{ccccccccc}
0 & \rightarrow & H_n(V'_t) & \rightarrow & \hat{H} &
\rightarrow & H_{n-1}(V_t) & \rightarrow & 0 \\
0 & \rightarrow & H_{n+1}(V^{(2)}_t) & \rightarrow &
\hat{H}' & \rightarrow & H_n(V'_t) & \rightarrow & 0 \\
 & &  \vdots & & \vdots & & \vdots & & \\
0 & \rightarrow & H_{n+p-1}(V^{(p)}_t) & \rightarrow &
\hat{H}^{(p-1)} &
\rightarrow & H_{n+p-2}(V^{(p-1)}_t) & \rightarrow & 0 \\
 & &  0 & \rightarrow & \hat{H}^{(p)} & \rightarrow
& H_{n+p-1}(V^{(p)}_t) & \rightarrow & 0
\end{array}
$$
They give rise to a long exact sequence (cf.\ \cite[p.~163]{AGV})
$$ 0 \rightarrow \hat{H}^{(p)} \rightarrow \hat{H}^{(p-1)} \rightarrow \ldots
\rightarrow \hat{H}' \rightarrow \hat{H} \rightarrow H_{n-1}(V_t)
\rightarrow 0 .$$
Each of the modules in this sequence is a free $\ZZ$-module of finite
rank. Let
$\nu_i:= {\rm rank}\, \hat{H}^{(i)}$. Then
$$\mu = {\rm rank}\, H_{n-1}(V_t) = \sum_{i=0}^{p} (-1)^i \nu_i.$$
On each of the modules we have an intersection form
$\langle\cdot,\,\cdot\rangle$ defined as
in \cite{Eb}. On the module $H_{n+i-1}(V^{(i)}_t)$ ($0\le i\le p$) it is
the usual intersection form; on the module $\hat{H}^{(i)}$ it is the
pullback of the intersection form by the natural (boundary) homomorphism
$\hat{H}^{(i)}\to H_{n+i-1}(V^{(i)}_t)$. The form on
$H_{n-1}(V_t)$ is symmetric if $n$ is odd and skew-symmetric if
$n$ is even. The form on $\hat{H}^{(i)}$ is symmetric if $n+i$ is odd and
skew-symmetric if $n+i$ is even.

Denote by $\hat{H}^\ast$ the dual module of $\hat{H}=H_n(V'_t,V_t)$ and
let $(\cdot,\,\cdot) : \hat{H}^\ast \times \hat{H} \to \ZZ$ be the
Kronecker pairing. We want to define a variation operator
${\rm Var}: \hat{H}^\ast \to
\hat{H}$ or rather its inverse $\Var : \hat{H} \to \hat{H}^\ast$. For this
purpose we need the notion of a distinguished basis of thimbles.

Let $\tilde{f}_{p+1}:V'_t\to \CC$ be a generic perturbation of the
restriction of the function $f_{p+1}$ to $V'_t$ which has only non-degenerate
critical points with different critical values $z_1, \ldots , z_\nu$
($\nu=\nu_0$). Let $z_0$ be a non-critical value of $\tilde{f}_{p+1}$ with
$\| z_0 \| > \| z_j\|$ for $j=1, \ldots , \nu$. The level set
$\{ x\in V'_t : \tilde{f}_{p+1}(x) = z_0 \}$ is diffeomorphic to the Milnor
fibre $V_t$ of the {\em icis} $(V,0)$. Let $u_j$, $j=1, \ldots , \nu$, be
non-self-intersecting paths connecting the critical values $z_j$ with the
non-critical value $z_0$ in such a way that they lie inside the disc $D_{\| z_0
\|} = \{ z \in \CC : \| z \| \leq \| z_0 \| \}$ and every two of them intersect
each other only at the point $z_0$. We suppose that the paths $u_j$ (and
correspondingly the critical values $z_j$) are numbered clockwise according to
the order in which they arrive at $z_0$ starting from the boundary of the disc
$D_{\| z_0 \|}$. Each path $u_j$ defines up to orientation a thimble
$\del_j$ in the relative homology group $\hat{H}$. The system $\{ \del_1 ,
\ldots , \del_\nu \}$ is a basis of $\hat{H}$. A basis obtained in this way is
called {\em distinguished}. The self-intersection number of a thimble $\del$ is
equal to
$$\langle \del, \del \rangle = (-1)^{n(n-1)/2}(1+(-1)^{n-1}).$$
The {\em Picard-Lefschetz transformation} $h_{\del}:\hat{H} \to \hat{H}$
corresponding to the thimble
$\del$ is given by (cf.\  \cite{Eb})
$$h_{\del} (y) = y + (-1)^{n(n+1)/2}\langle y,\del \rangle \del \quad {\rm for}
\ y \in \hat{H}.$$

Going once around the disc $D_{\| z_0 \|}$ in the positive direction
(counterclockwise) along the boundary induces an automorphism of $\hat{H}$,
the {\em (classical) monodromy operator} $h_\ast$. If $\{ \del_1,
\ldots, \del_\nu\}$ is a distinguished basis of $\hat{H}$, then the
monodromy operator is given by
$$h_\ast = h_{\del_1} \circ h_{\del_2} \circ \cdots \circ h_{\del_\nu}.$$

\begin{definition} Let $\{ \del_1, \ldots , \del_\nu \}$ be a
distinguished basis of thimbles of $\hat{H}$ and let $\{ \nabla_1, \ldots ,
\nabla_\nu \}$ be the corresponding dual basis of $\hat{H}^\ast$. The
linear operator $\Var : \hat{H} \to \hat{H}^\ast$ (inverse of the {\em
variation operator}) is defined by
$$ \Var(\del_i) = (-1)^{n(n+1)/2} \nabla_i - \sum_{j < i} \langle \del_i,
\del_j \rangle \nabla_j.$$
\end{definition}

\begin{proposition}
The definition of the operator $\Var$ does not depend on the choice of the
distinguished basis.
\end{proposition}

\begin{proof}
Any two distinguished bases of thimbles can be transformed into each
other by the braid group transformations $\alpha_j$, $j=1, \ldots ,
\nu-1$, and by changes of orientations (see, e.g., \cite{AGV},
\cite{Eb}). Here the operation $\alpha_j$ is defined as follows:
$$ \alpha_j (\del_1, \ldots , \del_\nu) = (\delp_1, \ldots ,
\delp_\nu)$$ where $\delp_j = h_{\del_j}(\del_{j+1})= \del_{j+1} +
(-1)^{n(n+1)/2}\langle \del_{j+1},\del_{j} \rangle \del_j$,
$\delp_{j+1} = \del_j$, and
$\delp_i = \del_i$ for $i \neq j, j+1$.

It is easily seen that the definition of $\Var$ is invariant under a change of
orientation. Therefore it suffices to show that the definition of $\Var$ is
invariant under the transformation $\alpha_j$. One easily computes:
\begin{eqnarray*}
\langle \delp_r ,\delp_s \rangle & = & \langle \del_r, \del_s \rangle \quad
\mbox{for} \ 1 \leq r,s \leq \nu, \ r,s \neq j,j+1, \\
\langle \delp_j ,\delp_{j+1} \rangle & = &  -\langle \del_j, \del_{j+1}
\rangle , \\
\langle \delp_r ,\delp_j \rangle & = &  \langle \del_r,
\del_{j+1} \rangle +(-1)^{n(n+1)/2}\langle \del_{j+1}, \del_j \rangle \langle
\del_r, \del_j \rangle \ \mbox{ for} \ r \neq j,j+1, \\
\langle \delp_r ,\delp_{j+1} \rangle & = &  \langle \del_r, \del_j \rangle
\quad \mbox{for} \ r \neq j,j+1.
\end{eqnarray*}

Let $(\nabla'_1, \ldots , \nabla'_\nu)$ be the dual basis corresponding
to $(\delp_1, \ldots ,\delp_\nu)$.
Then
\begin{eqnarray*}
\nabla_{j+1} & = & \nabla'_j \\
\nabla_{j} & = & \nabla'_{j+1} + (-1)^{n(n+1)/2} \langle \del_{j+1} ,
\del_j \rangle \nabla'_{j}.
\end{eqnarray*}
One has
\begin{eqnarray*}
\Var(\delp_j) & = & 
\Var(\del_{j+1}) + (-1)^{n(n+1)/2} \langle \del_{j+1} , \del_j
\rangle \Var(\del_j) \\
& = & (-1)^{n(n+1)/2} \nabla_{j+1} - \sum_{k<j+1} \langle \del_{j+1} ,
\del_k \rangle \nabla_k \\
& & + (-1)^{n(n+1)/2}\langle \del_{j+1} , \del_j \rangle
((-1)^{n(n+1)/2}\nabla_j - \sum_{k<j} \langle \del_j , \del_k \rangle
\nabla_k ) \\
& = & (-1)^{n(n+1)/2}\nabla_{j+1} \\
& &  - \sum_{k<j} (\langle \del_{j+1} ,
\del_k \rangle + (-1)^{n(n+1)/2} \langle \del_{j+1} , \del_j \rangle
\langle \del_j , \del_k \rangle ) \nabla_k \\
& = & (-1)^{n(n+1)/2}\nabla'_j - \sum_{k<j} \langle \delp_j ,\delp_k
\rangle \nabla'_k,
\end{eqnarray*}
\begin{eqnarray*}
\Var(\delp_{j+1}) & = & 
(-1)^{n(n+1)/2}\nabla_j - \sum_{k<j} \langle \del_j , \del_k
\rangle \nabla_k \\
& = & (-1)^{n(n+1)/2}(\nabla_j - (-1)^{n(n+1)/2}\langle \del_{j+1} ,
\del_j \rangle \nabla_{j+1} ) \\
& & + \langle \del_{j+1} ,
\del_j \rangle \nabla_{j+1} - \sum_{k<j} \langle \del_j , \del_k
\rangle \nabla_k \\
& = &  (-1)^{n(n+1)/2}\nabla'_{j+1} - \sum_{k<j+1} \langle \delp_{j+1}
, \delp_k \rangle \nabla'_k .
\end{eqnarray*}
For $i>j+1$ we have
\begin{eqnarray*}
\Var(\delp_i) & = &
(-1)^{n(n+1)/2}\nabla_i - \sum_{k<i} \langle \del_i , \del_k
\rangle \nabla_k \\
& = & (-1)^{n(n+1)/2}\nabla'_i - \hspace{-1.5mm}
\sum\limits_{\substack{k<i \\ k \neq j,j+1}} \langle \delp_i
, \delp_k \rangle \nabla'_k -\langle
\del_i , \del_{j+1} \rangle \nabla_{j+1} - \langle \del_i, \del_j
\rangle \nabla_j \\
& = & (-1)^{n(n+1)/2}\nabla'_i - \hspace{-1.5mm}
\sum\limits_{\substack{k<i \\ k \neq j,j+1}} \! \! \langle \delp_i ,
\delp_k
\rangle \nabla'_k \\ & & - \langle \delp_i, \delp_j \rangle \nabla'_j  +
(-1)^{n(n+1)/2}\langle \del_{j+1} , \del_j \rangle \langle \del_i ,
\del_j \rangle \nabla'_j \\
& & - \langle \delp_i , \delp_{j+1} \rangle
\nabla'_{j+1} - (-1)^{n(n+1)/2}\langle \del_{j+1} , \del_j \rangle
\langle \delp_i , \delp_{j+1} \rangle \nabla'_j \\
& = & (-1)^{n(n+1)/2}\nabla'_i - \sum_{k<i} \langle
\delp_i , \delp_k \rangle \nabla'_k
\end{eqnarray*}
The corresponding formula for $i<j$ is obvious.
\end{proof}

\begin{remark}
There is an interesting problem to give an invariant
(topological) definition of the variation operator.
\end{remark}

Let $S: \hat{H} \to \hat{H}^\ast$ be the mapping defined by the
intersection form on $\hat{H}$: $(Sx, y) =\langle x , y \rangle$,
$x,\,y\in \hat{H}$. The mapping $\Var$ is defined
in such a way that one has the equality $S = -\Var + (-1)^n (\Var)^T$
where $(\Var)^T$ means the transpose operator $(\Var)^T :
\hat{H}^{\ast\ast}= \hat{H} \to \hat{H}^\ast$.

\begin{remark}
We emphasize that the intersection number $\langle \del_i , \del_j
\rangle$ is the entry of the matrix of the operator $S$ with {\em
column} index $i$ and {\em row} index $j$: see the remark in
\cite[p.~45]{AGV}.
\end{remark}

\begin{remark}
The operator $V: \hat{H}\to \hat{H}^\ast$ defined in
\cite[p.~18]{Eb} differs from $\Var$ by sign.
\end{remark}

\begin{proposition}
The monodromy operator $h_\ast$ and the operator $\Var$ are related by the
following formula:
$$h_\ast = (-1)^n {\rm Var} (\Var)^T. $$
\end{proposition}

\begin{proof}
This can be computed directly; see also \cite[Chap.~V, \S 6,
Exercice 3]{B}.
\end{proof}

In the same way, for each $i$ with $1 \leq i \leq p$ an operator
${\rm Var}^{-1}_i : \hat{H}^{(i)} \to (\hat{H}^{(i)})^\ast$ is defined.

\section{The index of the gradient vector field on an isolated
complete intersection singularity}\label{sec3}

Let $(V', 0)=\{f_1=f_2=\ldots=f_p=0\}\subset(\CC^{n+p}, 0)$ be a real
$n$-dimensional
{\em icis} (it means that the function germs $f_i:(\CC^{n+p},
0)\to(\CC, 0)$ are real). We assume that its real part $V'\cap\RR^{n+p}$
does not coincide with the origin (and thus is $n$-dimensional). Let
$g=f_{p+1}:(\CC^{n+p}, 0)\to(\CC, 0)$ be a germ of a real analytic function
such that its restriction to $V'\setminus\{0\}$ has no critical points. A
Riemannian metric on $\RR^{n+p}$ determines the gradient vector field
$X=\mbox{grad}\,g$ of the restriction of the function $g$ to
$(V'\cap\RR^{n+p})\setminus\{0\}$. This vector field has no singular
points on a punctured neighbourhood of the origin in
$V'\cap\RR^{n+p}$. Since the space of Riemannian metrics is connected,
the index \,$\mbox{ind}_0\,X$ of the gradient vector field doesn't
depend  on the choice of a metric. In the case $p=0$ (and thus
$V'=\CC^{n}$) the index of the gradient vector field of a function
germ $g$ can be expressed in terms of the action of the complex
conjugation on the Milnor lattice of the singularity $g$ (\cite{GZ},
\cite{V}). We give a generalisation of such a formula for {\em icis}.

Let $0<\eps_1\ll\eps_2\ll\ldots\ll\eps_{p+1}$ be real and small enough, let
$s=(s_1,\,\ldots,\, s_{p+1})$ with $s_i=\pm 1$. For $0\le i\le p$, let
$\hat{H}^{(i)}=\hat{H}^{(i)}_{s\eps}$ be the corresponding space of thimbles:
$\hat{H}^{(i)}=H_{n+i}(V^{(i+1)}_{s\eps}, V^{(i)}_{s\eps})$ for $0\le i\le
p-1$ ($s\eps=(s_1\eps_1, s_2\eps_2, \ldots, s_{p+1}\eps_{p+1}))$; see
Section~\ref{sec2} for $i=p$. Let $\sigma^{(i)}_s$ be the action of the
complex conjugation on the space $\hat{H}^{(i)}$, let
$\Var_i:\hat{H}^{(i)}\to (\hat{H}^{(i)})^\ast$
be the inverse of the corresponding variation operator. The operator
$\Var_i\sigma^{(i)}_s$ acts from the space $\hat{H}^{(i)}$ to its dual
$(\hat{H}^{(i)})^\ast$ and thus defines a bilinear form on $\hat{H}^{(i)}$.

\begin{theorem}\label{theo1}
The bilinear forms $\Var_i\sigma^{(i)}_s$ are symmetric and
non-degener\-ate, and we have
\begin{eqnarray}\label{eq2}
{\rm ind}_0\,{\rm grad}\,g & = & s_{p+1}^{n}(-1)^{\frac{n(n+1)}{2}}
{\rm sgn}\,\Var_{}\sigma^{}_s \nonumber \\
 & & +\sum\limits_{i=1}^p(-1)^{\frac{(n+i)(n+i+1)}{2}}{\rm
sgn}\,\Var_{i}\sigma^{(i)}_s.
\end{eqnarray}
\end{theorem}

\begin{corollary}
The right-hand side of the equation {\rm (\ref{eq2})} does not depend on
$s=(s_1,\,\ldots,\, s_{p+1})$.
\end{corollary}

\begin{proof} Let us consider the restriction of the function $f_i$ to the
manifold $V^{(p-i+2)}_{s\eps}$. It may have degenerate critical points. Let
$\widetilde f_i:V^{(p-i+2)}_{s\eps}\to\CC$ be its real morsification (i.e., a
perturbation of $f_i$ which is a Morse function on $V^{(p-i+2)}_{s\eps}$ and
maps its real part $V^{(p-i+2)}_{s\eps}\cap\RR^{n+p}$ to $\RR$). For
$c\in\RR$, let $M_c^{(i)}=\{x\in
V^{(p-i+2)}_{s\eps}\cap\RR^{n+p}:\widetilde f_i\le c\}$. The topological
space $M_c^{(i)}$ is homotopy equivalent to
$V^{(p-i+2)}_{s\eps}\cap\RR^{n+p}$ or to
$V^{(p-i+1)}_{(s_1\eps_1,\ldots, s_{i-1}\eps_{i-1}, -\eps_i)}\cap\RR^{n+p}$
for $c$ greater than or less than all the critical values of $\widetilde f_i$
respectively. The standard arguments of Morse theory give
$$
\chi(V^{(p-i+2)}_{s\eps}\cap\RR^{n+p})=\chi(V^{(p-i+1)}_{(s_1\eps_1,\ldots,
s_{i-1}\eps_{i-1}, -\eps_i)}\cap\RR^{n+p}) + \sum\limits_{Q\in
S_{{\rm grad}\,\widetilde f_i}}{\rm ind}_Q\,{\rm grad}\,\widetilde f_i.
$$
Applying the same reasonings to the function $-\widetilde f_i$ one has
\begin{eqnarray*}
\chi(V^{(p-i+2)}_{s\eps}\cap\RR^{n+p}) & =
& \chi(V^{(p-i+1)}_{(s_1\eps_1,\ldots, s_{i-1}\eps_{i-1},
\eps_i)}\cap\RR^{n+p}) \\
 & & + (-1)^{n+p-i+1}\sum\limits_{Q\in
S_{{\rm grad}\,\widetilde f_i}}{\rm ind}_Q\,{\rm grad}\,\widetilde f_i.
\end{eqnarray*}
Thus
$$
\chi(V^{(p-i+2)}_{s\eps}\cap\RR^{n+p})=\chi(V^{(p-i+1)}_{s\eps}\cap\RR^{n+p}) +
(-s_i)^{n+p-i+1}
\sum\limits_{Q\in
S_{{\rm grad}\,\widetilde f_i}}{\rm ind}_Q\,{\rm grad}\,\widetilde f_i.
$$
(cf.\ \cite[Lemma in \S 2]{A}). From Proposition 1 one has
\begin{eqnarray*}
\lefteqn{{\rm ind}_{0}{\rm grad}\, g} \\
& & =\hspace{-3mm}\sum\limits_{Q\in
S_{{\rm grad}\,\widetilde f_{p+1}}}\hspace{-3mm}{\rm ind}_Q{\rm
grad}\,\widetilde
f_{p+1} - \chi(V'_{s\eps}\cap\RR^{n+p}) + 1  \\
& & =\hspace{-3mm}\sum\limits_{ Q\in S_{{\rm grad}\,\widetilde
f_{p+1}}}\hspace{-4.5mm}{\rm ind}_Q{\rm grad}\,\widetilde f_{p+1}
+ (-s_p)^{(n+1)}\hspace{-3mm}\sum\limits_{Q\in S_{{\rm grad}\,\widetilde
f_{p}}}\hspace{-3.5mm}{\rm ind}_Q{\rm grad}\,\widetilde f_{p}
- \chi(V''_{s\eps}\cap\RR^{n+p}) + 1 \\
& & = \ldots =\\
& & =  \hspace{-3mm}\sum\limits_{Q\in S_{{\rm grad}\,\widetilde
f_{p+1}}} \hspace{-3mm}{\rm ind}_Q{\rm grad}\,\widetilde f_{p+1}
+ \sum\limits_{i=1}^p(-s_{p-i+1})^{(n+i)}\hspace{-3mm}\sum\limits_{Q\in
S_{{\rm grad}\,\widetilde f_{p-i+1}}}\hspace{-3mm}{\rm ind}_Q{\rm
grad}\,\widetilde
f_{p-i+1}.
\end{eqnarray*}

Now Theorem~\ref{theo1} follows from the following statement.
\end{proof}

\begin{theorem}\label{theo2}
$$
\sum\limits_{Q\in S_{{\rm grad}\,\widetilde f_{p-i+1}}}
{\rm ind}_Q{\rm grad}\,\widetilde f_{p-i+1}
=(s_{p-i+1})^{n+i}(-1)^{\frac{(n+i)(n+i+1)}{2}}{\rm
sgn}\,\Var_i\sigma^{(i)}_s.
$$
\end{theorem}

\begin{proof}
Let us suppose that $s_{p-i+1} =1$. The case $s_{p-i+1}=-1$ can be
reduced to this by multiplying the function $f_{p-i+1}$ (and the function
$\widetilde{f}_{p-i+1}$) by $(-1)$. Let $\widetilde{s} = (s_1, \ldots ,
s_{p-i}, -1, s_{p-i+2}, \ldots, s_{p+1})$. Without loss of generality we
can suppose that all  critical values of the function
$\widetilde{f}_{p-i+1}$ lie inside the circle
$\{z: \| z \|\le\frac{\eps_{p-i+1}}{2} \}$ and have different real parts
(except, of course, values at complex conjugate points). Let us identify
the space
$\hat{H}^{(i)}_{\widetilde{s}}$ with the space $\hat{H}^{(i)}_s$ using a
path which connects $-\eps_{p-i+1}$ with $+\eps_{p-i+1}$ in the upper
half plane outside the circle
$\{z: \| z\| < \frac{\eps_{p-i+1}}{2} \}$ (e.g., the half circle
$\{z: \| z\| = \eps_{p-i+1}\}$). This identification permits to consider
$\sigma_s^{(i)}$ and
$\sigma_{\widetilde{s}}^{(i)}$ as operators on the space $\hat{H}^{(i)}
=\hat{H}^{(i)}_{s\eps}$. Just as in \cite{GZ} the classical monodromy
operator $h_{\ast}^{(i)}: \hat{H}^{(i)} \to \hat{H}^{(i)}$ can be
represented in the form
\begin{equation}\label{eq3}
h_{\ast}^{(i)} = \sigma^{(i)}_s \sigma^{(i)}_{\widetilde{s}}.
\end{equation}
A distinguished basis of the space $\hat{H}^{(i)}$ is defined by a
system of paths connecting the critical values of the function
$\widetilde{f}_{p-i+1}$ with the non-critical value $\eps_{p-i+1}$. Let
us choose the following system of paths (cf.\ Fig.~\ref{fig1}).
\begin{figure}
\vspace{2cm}
\centering
\unitlength1cm
\begin{picture}(9,5)
\put(0,0){\includegraphics{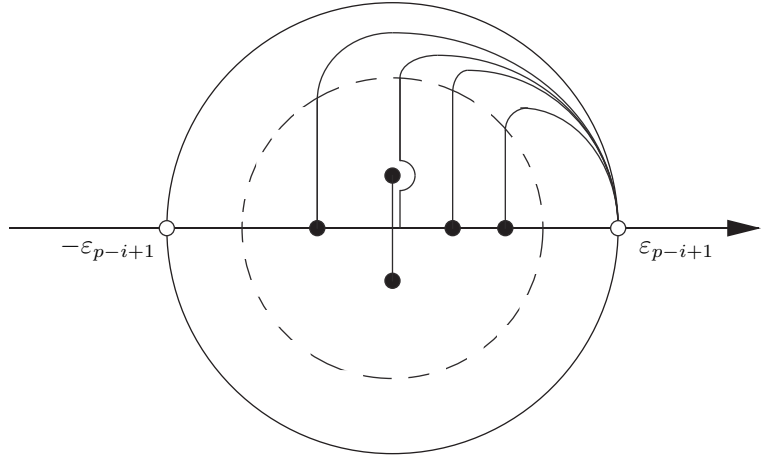}}
\put(0.7,2.7){$-\eps_{p-i+1}$}
\put(8.4,2.7){$\eps_{p-i+1}$}
\end{picture}
\caption{The choice of paths}\label{fig1}
\end{figure}
The paths from real critical values go vertically upwards
up to the boundary of the circle $\{ z: \| z \| \leq
\frac{\eps_{p-i+1}}{2} \}$. The paths from complex conjugate critical
values go vertically (upwards or downwards to the real axis and then go
vertically upwards to the boundary of the circle $\{ z: \| z \| \leq
\frac{\eps_{p-i+1}}{2} \}$ avoiding from the right side a neighbourhood
of the critical value with positive imaginary part. From the boundary of
the circle $\{ z: \| z \| \leq \frac{\eps_{p-i+1}}{2} \}$ all the paths
go to the non-critical value $\eps_{p-i+1}$ in the upper half plane (see
Fig.~\ref{fig1}). The cycles are ordered in the usual way which in this
case means that they follow each other in the order of decreasing real
parts of the corresponding critical values; the vanishing cycle
corresponding to the critical value with negative imaginary part
precedes that with the positive one.

In the sequel we shall consider the matrices of the operators
$\sigma^{(i)}_s$, $\sigma^{(i)}_{\widetilde{s}}$, $\Var_i$, etc.\  as block
matrices with blocks of size $1 \times 1$, $1 \times 2$, $2 \times 1$, and $2
\times 2$ corresponding to real critical values and to pairs of complex
conjugate critical values of the function $\widetilde{f}_{p-i+1}$. The
matrix of the operator $\sigma_s^{(i)}$ is an upper triangular block
matrix. Its diagonal entry corresponding to a real critical value is
equal to
$(-1)^m$ where $m$ is the Morse index of the critical point. A diagonal
block of size $2 \times 2$ corresponding to a pair of complex conjugate
critical values is equal to
$\left( \begin{array}{cc} 0 & 1 \\ 1 & 0 \end{array} \right)$. The matrix of
the operator
$\sigma_{\widetilde{s}}^{(i)}$ is a lower triangular block matrix (we do
not need a precise description of its diagonal blocks). The matrix of the
operator
$\Var_i$ is upper triangular with diagonal entries equal to
$(-1)^{(n+i)(n+i+1)/2}$
(the dual of the space $\hat{H}^{(i)}$ is
endowed with the basis dual to the one of $\hat{H}^{(i)}$). One has
$h_{\ast}^{(i)} = \sigma_s^{(i)} \sigma_{\widetilde{s}}^{(i)} =
(-1)^{n+i} {\rm Var}_i (\Var_i)^T$. Thus
$\Var_i \sigma_s^{(i)} = (-1)^{n+i} (\Var_i)^T \sigma_{\widetilde{s}}^{(i)}$.
The matrices $\Var_i\sigma_s^{(i)}$ and $(\Var_i)^T
\sigma_{\widetilde{s}}^{(i)}$
are upper triangular and lower triangular respectively. Thus the
matrix $\Var_i \sigma_s^{(i)}$ is in fact block diagonal with
the diagonal entry
$(-1)^{((n+i)(n+i+1)/2) +m}$
corresponding to a real
critical point of the function $\widetilde{f}_{p-i+1}$ ($m$ is the Morse
index) and with the diagonal block of the form
$$(-1)^{\frac{(n+i)(n+i+1)}{2}} \left( \begin{array}{cc} a & 1 \\ 1 & 0
\end{array} \right) $$
corresponding to a pair of complex conjugate critical points (up to a sign
$a$ is the intersection number of the corresponding cycles). This description
implies Theorem~\ref{theo2}.
\end{proof}

The formula~(\ref{eq2}) expresses the index of a gradient vector
field in terms of bilinear forms on the spaces of thimbles. Actually
each second summand of it can be expressed in terms of bilinear forms on
the corresponding spaces of vanishing cycles. Let $\Sigma_s^{(i)}$ be the
quadratic form on the space $H_{n+i-1}(V_{s\eps}^{(i)})$ of vanishing
cycles defined by
$\Sigma_s^{(i)}(x,y)=\langle\sigma_s^{(i)}x,y\rangle$. As above, let
$\widetilde{s} = (s_1,\linebreak[0] \ldots,\linebreak[0]
s_{p-i},\linebreak[0] -1,\linebreak[0] s_{p-i+2},
\ldots, s_{p+1})$.

\begin{theorem}\label{theo3} For $n+i$ odd
$$
\sum\limits_{Q\in S_{{\rm grad}\,\widetilde f_{p-i+1}}}
{\rm ind}_Q{\rm grad}\,\widetilde f_{p-i+1}
=s_{p-i+1}(-1)^{\frac{n+i+1}{2}}({\rm
sgn}\,\Sigma_{\widetilde s}^{(i)}-{\rm sgn}\,\Sigma_s^{(i)})/2.
$$
\end{theorem}

\begin{proof}
It is essentially the same as in \cite{GZ}. One has to notice that the
kernel of the natural (boundary) homomorphism
$\hat{H}^{(i)}_{s\eps}\to H_{n+i-1}(V_{s\eps}^{(i)})$ is contained in the
kernel of the quadratic form $\Sigma_s^{(i)}$ and thus
\ ${\rm sgn}\,\Sigma_{\widetilde s}^{(i)}$ coincides with the signature
of the form $\langle\sigma_s^{(i)}\cdot,\cdot\rangle$ on the space
$\hat{H}^{(i)}_{s\eps}$.
\end{proof}

\begin{remark} For $n+i$ even, it is not possible to express the number
$$\sum\limits_{Q\in S_{{\rm grad}\,\widetilde f_{p-i+1}}}
{\rm ind}_Q{\rm grad}\,\widetilde f_{p-i+1}$$ in terms of invariants
defined by the space of vanishing cycles. It can be understood from the
following example. Let $n=2$, $p=1$, $f_1=x_1^2+x_2^2-x_3^2$, $f_2=x_3$,
$i=0$. The discussed sum is different for $s_1=1$ and for $s_1=-1$
(i.e., for $t_1=s_1\eps_1$ positive or negative). However the line
$\ell=\{t_1=0\}$ is not in the bifurcation set for vanishing cycles: it
doesn't lie in the discriminant of the map $(f_1, f_2):(\CC^3,
0)\to(\CC^2, 0)$. On the other hand the discriminant of the map
$f_1:(\CC^3, 0)\to(\CC, 0)$ coincides with the origin $0\in\CC$ and
thus the line $\ell$ is in the bifurcation set for thimbles.
\end{remark}

\begin{remark}
It seems that a relation between complex conjugation
and the monodromy operator (similar to (\ref{eq3})) was first used in
\cite{A'C}. In \cite{GZ} it was written in an explicit way. In
\cite{MP} a similar relation was used to find some properties of
the Euler characteristics of links of complete intersection varieties.
This paper partially inspired our work.
\end{remark}

\end{document}